\documentclass[11pt]{article}
\usepackage[T1]{fontenc}
\usepackage[utf8]{inputenc}
\usepackage{lmodern,microtype}
\usepackage[a4paper,margin=30mm]{geometry}
\usepackage{amsmath,amssymb,amsthm,mathtools,mathrsfs,array,booktabs,needspace}
\usepackage{aliascnt}
\usepackage{xcolor,xurl,hyperref}
\usepackage[nameinlink,capitalise,noabbrev]{cleveref}
\usepackage{fancyhdr}
\hypersetup{colorlinks=true,linkcolor=blue!45!black,citecolor=green!35!black,
 urlcolor=blue!55!black,pdftitle={Historical Gromov--Witten classes and Kodaira dimension},
 pdfsubject={Symplectic deformation invariance}}
\newtheorem{theorem}{Theorem}[section]
\newaliascnt{question}{theorem}\newtheorem{question}[question]{Question}\aliascntresetthe{question}
\crefname{question}{Question}{Questions}\Crefname{question}{Question}{Questions}
\crefname{theorem}{Theorem}{Theorems}\Crefname{theorem}{Theorem}{Theorems}
\crefname{proposition}{Proposition}{Propositions}\Crefname{proposition}{Proposition}{Propositions}
\crefname{lemma}{Lemma}{Lemmas}\Crefname{lemma}{Lemma}{Lemmas}
\crefname{corollary}{Corollary}{Corollaries}\Crefname{corollary}{Corollary}{Corollaries}
\crefname{remark}{Remark}{Remarks}\Crefname{remark}{Remark}{Remarks}
\crefname{section}{Section}{Sections}\Crefname{section}{Section}{Sections}
\crefname{equation}{equation}{equations}\Crefname{equation}{Equation}{Equations}
\newaliascnt{proposition}{theorem}\newtheorem{proposition}[proposition]{Proposition}\aliascntresetthe{proposition}
\newaliascnt{lemma}{theorem}\newtheorem{lemma}[lemma]{Lemma}\aliascntresetthe{lemma}
\newaliascnt{corollary}{theorem}\newtheorem{corollary}[corollary]{Corollary}\aliascntresetthe{corollary}
\newaliascnt{conjecture}{theorem}\newtheorem{conjecture}[conjecture]{Conjecture}\aliascntresetthe{conjecture}
\crefname{conjecture}{Conjecture}{Conjectures}\Crefname{conjecture}{Conjecture}{Conjectures}
\theoremstyle{remark}
\newaliascnt{remark}{theorem}\newtheorem{remark}[remark]{Remark}\aliascntresetthe{remark}
\newcommand{\Q}{\mathbb Q}
\newcommand{\PP}{\mathbb P}\newcommand{\cO}{\mathcal O}

\newcommand{\Exc}{\operatorname{Exc}}\newcommand{\Supp}{\operatorname{Supp}}

\title{\bfseries On the symplectic deformation invariance of Kodaira dimension for threefolds}
\author{Jingcao Wu\thanks{School of Mathematics, Shanghai University of Finance and Economics, Shanghai 200433, People's Republic of China. E-mail: \texttt{wujincao@shufe.edu.cn}.}}
\date{}
\begin{document}\maketitle
\begin{abstract}
We study the symplectic deformation invariance of Kodaira dimension for smooth projective threefolds. We give a cohomological criterion for recovering the pullback of the canonical class of a smooth model and apply it to iterated point blow-ups. This yields an invariance result for non-uniruled threefolds whose canonical classes are divisible by two modulo torsion. We also discuss local Gromov--Witten invariants of divisorial contractions and the additional questions that arise for curve centres.
\end{abstract}

\section{Introduction}\label{sec:introduction}
In his survey on rational curves on complex manifolds, Cascini discusses which algebro-geometric properties are preserved under symplectic deformation \cite[Section~3.3]{Cascini13}. Uniruledness is a fundamental example: it is invariant under symplectic equivalence \cite[Theorem~3.6]{Cascini13}. Cascini also records the corresponding question for rational connectedness as Conjecture~3.8. It is known in dimensions two and three: the surface case follows from the classification of surfaces and the uniruledness theorem, while the threefold case was established by Tian, building on partial results of Voisin \cite[Conjecture~3.8]{Cascini13} \cite{Voisin08,Tian12}. These results motivate the following question about the Kodaira dimension itself.
\begin{question}\label{que:main}
Let $X$ and $Y$ be smooth complex projective varieties whose K\"ahler symplectic structures are deformation equivalent. Is $\kappa(X)=\kappa(Y)$?
\end{question}
The answer is already known in complex dimensions one and two. For curves it is determined by the genus. For surfaces, Tian-Jun Li defined the symplectic Kodaira dimension of a closed symplectic four-manifold and proved that it depends only on the oriented diffeomorphism type; for a K\"ahler surface it agrees with the usual holomorphic Kodaira dimension \cite{Li06}. Dorfmeister and Zhang further established the agreement between the holomorphic and symplectic Kodaira dimensions for smooth four-manifolds carrying both complex and symplectic structures \cite{DZ09}. Since symplectic deformation equivalence gives an orientation-preserving diffeomorphism, these results settle \cref{que:main} for surfaces.

For smooth projective threefolds, the threefold minimal model program and the abundance theorem imply that a pseudo-effective canonical divisor has nonnegative Kodaira dimension: one passes to a minimal model with nef canonical divisor, which is semiample by the threefold abundance theorem \cite{Mori88,Kawamata92}. By the BDPP theorem, a smooth projective variety is uniruled if and only if its canonical divisor is not pseudo-effective \cite{BDPP13}. Thus, in dimension three, $\kappa=-\infty$ if and only if the variety is uniruled. The symplectic invariance of uniruledness \cite{Kollar98,Ruan99} therefore settles the case $\kappa=-\infty$ for smooth projective threefolds. Tian's theorem on the symplectic invariance of rational connectedness in dimension three gives a stronger statement within this already settled uniruled range \cite{Tian12}; it does not address the non-uniruled case. The remaining question in dimension three concerns the invariance of the actual nonnegative Kodaira dimension. For non-uniruled varieties, the difficulty is that the canonical class of a nonminimal model contains exceptional contributions, so its powers need not determine the Kodaira dimension.

Throughout, symplectic deformation equivalence includes a diffeomorphism identifying the deformation classes of the symplectic forms. Under this identification, the canonical classes, cup products, and genus-zero Gromov--Witten invariants agree. Unless otherwise specified, all cohomology and homology groups have rational coefficients.

\begin{theorem}\label{thm:main}
Let $X$ and $Y$ be symplectically deformation equivalent smooth projective threefolds.
\begin{enumerate}
\item If both are obtained from smooth projective varieties with nef canonical divisor by finite, possibly empty, sequences of point blow-ups, then $\kappa(X)=\kappa(Y)$.
\item If $X$ and $Y$ are non-uniruled and their canonical classes are divisible by two modulo torsion in integral cohomology, then $\kappa(X)=\kappa(Y)$.
\end{enumerate}
\end{theorem}
The evenness assumption in the second assertion of \cref{thm:main} excludes contraction types E1, E3, E4, and E5, since their primitive contracted rational curves have anticanonical degree one. Thus only a contraction of type E2, the blow-down of a smooth point, can occur in the reduction to a nef canonical divisor. In particular, the proof of the second assertion does not use the unresolved E4 virtual-cycle comparison, which remains relevant to the general negative-ray detection criterion in \cref{lem:nef}. The first assertion, on the other hand, does not follow from birational invariance alone, since the two nef models need not be birational.

The paper is organized as follows. The reconstruction criterion, the point-blow-up detection result, and the proof of \cref{thm:main} are given in \cref{sec:reconstruction}. The local calculations for divisorial contractions are collected in \cref{sec:local}. Historical classes and the remaining comparisons for curve centres and singular minimal models are discussed in \cref{sec:historical}.

\section{Negative Gromov--Witten classes, reconstruction, and point blow-ups}\label{sec:reconstruction}
Write $k_X=c_1(K_X)$, and let $V_X^-$ be the rational span in $H_2(X,\Q)$ of the nonzero integral classes $\beta$ such that $k_X\cdot\beta<0$ and some primary genus-zero Gromov--Witten invariant in class $\beta$ is nonzero. We allow any number of marked points and rational cohomology insertions. For a class considered modulo torsion, we sum the invariants over its finitely many torsion lifts. This convention is invariant under symplectic deformation and avoids identifying distinct integral classes solely on the basis of numerical equivalence.

\begin{lemma}\label{lem:nef}
If $K_X$ is nef, then $V_X^-=0$. Conversely, suppose that $X$ is a smooth non-uniruled projective threefold and that every negative divisorial extremal ray supports a nonzero primary genus-zero Gromov--Witten invariant in its contracted curve class. Then $V_X^-=0$ implies that $K_X$ is nef.
\end{lemma}
\begin{proof}
The first assertion follows because a class represented by a stable map has nonnegative intersection with a nef divisor. For the converse, suppose that $K_X$ is not nef. The cone and contraction theorems then give a $K_X$-negative extremal contraction. Since $X$ is non-uniruled, the contraction is not of fibre type; by the classification of contractions from smooth threefolds, it is one of the five divisorial types E1--E5. The detection hypothesis then places its contracted curve class in $V_X^-$, a contradiction.
\end{proof}
The local calculations in \cref{sec:local} establish this activity for E1, E2, E3, and E5. The E4 case requires the obstruction-theoretic comparison discussed there.

\begin{proposition}\label{prop:ring}
Let $\pi:X\to Z$ be a composition of blow-ups of smooth projective varieties of dimension $n\ge2$ along smooth connected centres of codimension at least two. Suppose $V_X^-=\ker(\pi_*:H_2(X)\to H_2(Z))$. Put
\[
U=(V_X^-)^\perp\subset H^2(X),\qquad W=\{v\in H^2(X):\int_X v u_1\cdots u_{n-1}=0\ \text{for all }u_i\in U\}.
\]
There is a unique $P\in U$ such that $P-k_X\in W$, and $P=\pi^*k_Z$.
\end{proposition}
\begin{proof}
The blow-up formula gives $H^2(X)=\pi^*H^2(Z)\oplus W_{\rm exc}$, where $W_{\rm exc}$ is spanned by the total transforms of the exceptional divisors. Duality identifies $U$ with $\pi^*H^2(Z)$. Since each exceptional divisor pushes forward to zero as a codimension-one cycle, the projection formula gives $W_{\rm exc}\subseteq W$. Conversely, write $v=\pi^*a+e$, with $e\in W_{\rm exc}$. If $v\in W$, then $\int_Z a b h^{n-2}=0$ for every $b\in H^2(Z)$ and every ample class $h$. The Hard Lefschetz theorem and Poincar\'e duality imply that $a=0$. Hence $W=W_{\rm exc}$, and the canonical divisor formula for blow-ups gives the result.
\end{proof}

\begin{corollary}\label{cor:reconstruct}
If $X$ and $Y$ are symplectically deformation equivalent smooth projective threefolds admitting morphisms as in \cref{prop:ring} to smooth varieties with nef canonical divisors, and both satisfy its kernel hypothesis, then their Kodaira dimensions coincide.
\end{corollary}
\begin{proof}
The reconstructed class $P$ is preserved by the cohomology identification. By the threefold abundance theorem \cite{Kawamata92}, the Kodaira dimension of a nef canonical divisor equals its numerical dimension, which is the largest $j$ for which $P^j\ne0$ (with $j=0$ allowed). The result follows from the birational invariance of Kodaira dimension.
\end{proof}

\begin{proposition}\label{prop:isolation}
Suppose $\pi:X\to Z$ is as in \cref{prop:ring}, $K_Z$ is nef, and no rational curve in $Z$ meets $\pi(\Exc(\pi))$. Then $V_X^-\subseteq\ker\pi_*$.
\end{proposition}
\begin{proof}
Write $K_X=\pi^*K_Z+E$, where $E$ is effective and exceptional. A stable map of negative canonical degree has a component on which $K_X$ has negative degree. Since $\pi^*K_Z$ is nef, this component must lie in $\Supp E$, so its image in $Z$ is a point of $\pi(\Exc(\pi))$. If another component had nonconstant image in $Z$, the first such component along a path in the connected dual tree would meet $\pi(\Exc(\pi))$, contrary to the hypothesis. Thus every component is contracted by $\pi$.
\end{proof}

\subsection{Point blow-ups and exceptional detection}
\begin{lemma}\label{lem:points}
Any ordered sequence of point blow-ups of a smooth connected projective variety can be deformed, with its historical exceptional classes marked, to a sequence of blow-ups at distinct general points of the original variety.
\end{lemma}
\begin{proof}
Start with the original variety as the fibre over a point. At each stage, take the total space of the current smooth projective family as the next parameter space, pull back the family, and blow up its diagonal section. Inductively, this constructs a smooth projective family over a smooth irreducible parameter space containing every ordered configuration. Configurations of distinct points form a dense open subset, and the total transforms of the exceptional divisors provide global markings.
\end{proof}
\begin{proposition}\label{prop:point}
Let $Z$ be a smooth, non-uniruled projective variety with nef canonical divisor, and let $X$ be obtained from $Z$ by a sequence of point blow-ups. Then $V_X^-=\ker\pi_*$.
\end{proposition}
\begin{proof}
By \cref{lem:points} and deformation invariance, it is enough to consider blow-ups at distinct very general points. Since $Z$ is non-uniruled, a very general point lies on no rational curve: each of the countably many components of the relevant spaces of rational maps has nondominant evaluation map. Thus \cref{prop:isolation} applies. Let $E\simeq\PP^{n-1}$ be an exceptional divisor and let $f$ be the class of a line in $E$. Pullback of an ample divisor on $Z$ shows that every stable map in class $f$ is contained in $E$. The space of lines in $E$ is the Grassmannian $\operatorname{Gr}(2,n)$. For a line $f\subset E$, the normal bundle sequence gives
\[
N_{f/X}\simeq \cO_{\PP^1}(1)^{\oplus(n-2)}\oplus\cO_{\PP^1}(-1),
\]
so $H^1(f,N_{f/X})=0$. The two-point moduli space is therefore unobstructed. Taking insertions of codimension $n-1$ that restrict to point classes on $E$, the corresponding invariant is one. Such insertions can be represented by suitable signed powers of the class of $E$. Hence every historical exceptional direction lies in $V_X^-$, and the reverse inclusion follows from \cref{prop:isolation}.
\end{proof}

\begin{proof}[Proof of \cref{thm:main}]
The first assertion follows from \cref{prop:point,cor:reconstruct}. For the second, let $V$ denote either $X$ or $Y$. If $K_V$ is not nef, the cone and contraction theorems give a $K_V$-negative extremal contraction. Since $V$ is non-uniruled, this contraction is divisorial. By the classification of negative extremal contractions from smooth threefolds \cite{Mori82}, it is of type E1--E5. In types E1, E3, E4, and E5, the primitive contracted rational curve $f$ satisfies $-K_V\cdot f=1$. This is impossible because $k_V$ is divisible by two modulo torsion in integral cohomology. Hence the contraction is of type E2: it is the blow-down $q:V\to V'$ of a smooth point, and
\[
K_V=q^*K_{V'}+2E.
\]
The integral blow-up decomposition of $H^2$ shows that divisibility of $k_V$ by two modulo torsion descends to $k_{V'}$. The contraction also preserves non-uniruledness. Since each contraction lowers the Picard number, iteration terminates at a smooth projective variety with nef canonical divisor. Thus both $X$ and $Y$ arise from smooth varieties with nef canonical divisor by finite sequences of point blow-ups, and assertion (1) applies. The unresolved E4 virtual-cycle comparison plays no role in this argument.
\end{proof}

\section{Local invariants of divisorial contractions}\label{sec:local}
We use the standard genus-zero virtual class and evaluation morphisms $\mathrm{ev}_i$ \cite{BF97}. In each calculation, the pullback of an ample divisor from the contraction target confines maps in a contracted class to the exceptional locus. Let $E$ denote the exceptional divisor and let $f$ be a minimal contracted curve. For types E1, E3, E4 and E5 define
\[
\Theta_f=\mathrm{PD}\, (\mathrm{ev}_1)_*[\overline{\mathcal M}_{0,1}(X,f)]^{\rm vir}\in H^2(X).
\]
For E2, use two markings and insert $u^2$ at the second marking, where $u$ is ample; normalize by dividing by $(u\cdot f)^2$. This also defines a divisor class $\Theta_f$.
\begin{proposition}\label{prop:local}
For types E1, E2, E3, and E5, the following local formulas hold, with torsion summation as above. The E4 entry is the expected value in the standard local model; it is valid only if the stable-map obstruction theory is identified with the zero-section obstruction theory described below.
\[
\begin{array}{c|ccccc}
 &\mathrm{E1}&\mathrm{E2}&\mathrm{E3}&\mathrm{E4}&\mathrm{E5}\\\hline
-K_X\cdot f&1&2&1&1&1\\
E\cdot f&-1&-1&-1&-1&-2\\
\Theta_f/[E]&1&1&2&2&-1
\end{array}
\]
For every type for which the indicated virtual-cycle calculation is established and $\Theta_f\cdot f\ne0$, the canonical pullback is
\[
q^*k_Z=k_X-\frac{k_X\cdot f}{\Theta_f\cdot f}\,\Theta_f.
\]
In particular, this formula is unconditional for E1, E2, E3 and E5, and conditional on the stated comparison for E4.
\end{proposition}
\begin{proof}
For E1, the exceptional divisor is a ruled surface over the smooth centre. A fibre has normal bundle $\cO\oplus\cO(-1)$, and the one-point evaluation pushes forward to the exceptional divisor. For E2, the moduli space of lines in $\PP^2$ is unobstructed because the normal bundle of the exceptional divisor is $\cO(-1)$. Two distinct points determine a unique line, so the two-point evaluation pushes forward to $[E\times E]$. The normalized insertion at the second marking then gives $\Theta_f=[E]$.

For E3, the exceptional divisor is a smooth quadric with normal bundle $\cO(-1,-1)$. Each ruling contributes one copy of $[E]$. The two ruling classes agree in rational homology: they lie on the same extremal numerical ray and have the same degree. On a smooth projective variety, numerical and homological equivalence agree for algebraic curve classes with rational coefficients, by the Hard Lefschetz theorem and the Lefschetz $(1,1)$ theorem.

For E5, the exceptional divisor is $\PP^2$ with normal bundle $\cO(-2)$. The unmarked moduli space is the dual plane $B$. Its obstruction line is $R^1p_*\cO_{\PP_B(S)}(-2)=\det S=\cO_B(-1)$, where $S$ is the tautological rank-two bundle and $p$ is the universal line. The resulting virtual class is $-H_B\cap[B]$, and the incidence correspondence pushes the marked virtual class forward to $-[E]$.

For E4, consider the blow-up of the origin in the hypersurface $xz-y^2+w^3=0$. In the ambient blow-up, the strict transform is given by $XZ-Y^2+rW^3=0$, where $r$ is the radial coordinate. The expected degree-one parameter scheme is the zero locus of the quadratic section of $\operatorname{Sym}^2S^*$ on $\operatorname{Gr}(2,4)$. In the chart $X=As+Bt$, $Y=Cs+Dt$, $Z=s$, and $W=t$, the equations are $A=C^2$, $B=2CD$, and $D^2=0$. Thus the chart has a doubled one-dimensional structure. This scheme-theoretic calculation does not, by itself, identify the stable-map virtual class. That requires a comparison of perfect obstruction theories, including the morphism to the cotangent complex, compatible with the local stable-map compactification and evaluation map. If such a comparison can be established, the local zero-section calculation predicts virtual multiplicity two and generic degree one for the one-point evaluation. The same issue arises for the standard E4 model with higher radial terms.

To obtain the canonical pullback formula, write $k_X=q^*k_Z+a[E]$. Pairing with $f$ gives $a=(k_X\cdot f)/(E\cdot f)$. If $\Theta_f=c[E]$ with $c\ne0$, then $\Theta_f\cdot f=c(E\cdot f)$, and therefore
\[
\frac{k_X\cdot f}{\Theta_f\cdot f}\Theta_f=\frac{k_X\cdot f}{E\cdot f}[E]=a[E].
\]
This proves the canonical pullback formula.
\end{proof}

\section{Historical classes and blow-up-centre comparisons}\label{sec:historical}
We begin with preservation and detection for point centres. We then formulate a reduced surface target and describe the additional obstruction-theoretic comparisons needed for curve centres. The conditional statements concerning curve centres are not used in the proof of \cref{thm:main}.
\subsection{Point-centre preservation and independent historical directions}
For a smooth blow-up $q:Y\to X$, write $q^!\alpha$ for the unique homological lift of $\alpha$ that pairs trivially with the exceptional divisor. For a sequence of blow-ups, use the composite lift. This notation is distinct from virtual pullback.
\begin{proposition}\label{prop:historical-point}
Suppose that every stable map in a nonzero genus-zero class $\alpha$ on a smooth projective variety $X$ has image contained in a proper closed subset $A$. Then, for any sequence of point blow-ups $\pi:Y\to X$, the primary invariants in class $\pi^!\alpha$ with insertions pulled back from $X$ agree with the corresponding invariants in class $\alpha$.
\end{proposition}
\begin{proof}
By \cref{lem:points}, we may deform the point configuration to distinct points outside $A$. Every stable map in the lifted class pushes forward, after stabilization, to a stable map in class $\alpha$, and its nonconstant components have image in $A$. Connectedness therefore rules out exceptional components over the chosen blow-up points. The two stable-map stacks are consequently isomorphic, including over nonreduced bases, because their maps lie in the common open set on which the blow-up is an isomorphism. Their standard obstruction theories agree on this open set. Deformation invariance gives the result for the original configuration.
\end{proof}
If the E4 local virtual-cycle calculation is established, \cref{prop:historical-point} also gives preservation of this contribution under a blow-up of the vertex, without a separate analysis of the new nonreduced moduli space. The historical class is then $h+f$, where $h$ is the strict ruling and $f$ is a line in the new exceptional plane. Its reduced moduli space parametrizes a point on the conic of tangent directions together with a line in the exceptional plane through that point. It is a $\PP^1$-bundle over the conic, hence isomorphic to $\PP^1\times\PP^1$: the quotient of $\cO^3$ by the tautological $\cO(-2)$ on the conic is $\cO(1)^{\oplus2}$. Under the assumed local virtual-cycle identification, the virtual pushforward to the original ruling parameter has coefficient two. A divisor insertion of degree one on the ruling then recovers the corresponding numerical contribution.

\begin{proposition}\label{prop:historical-directions}
Let $X_r\to\cdots\to X_0$ be a sequence of smooth blow-downs of points or curves, and let $f_i$ be the fibre class of $X_i\to X_{i-1}$. Suppose that, for each $i$, there is a class $\gamma_i$ on $X_r$ with negative canonical degree and a nonzero primary genus-zero Gromov--Witten invariant, whose pushforward to $X_i$ is a nonzero multiple of $f_i$. Then $\ker(H_2(X_r)\to H_2(X_0))\subset V_{X_r}^-$.
\end{proposition}
\begin{proof}
For $0\leq i\leq r$, set $K_i=\ker(H_2(X_r,\Q)\to H_2(X_i,\Q))$. Then $K_0=\ker(H_2(X_r)\to H_2(X_0))$ and $K_r=0$. By the blow-up formula, each quotient $K_{i-1}/K_i$ is one-dimensional, generated by the total-transform lift to $X_r$ of the exceptional fibre of $X_i\to X_{i-1}$. By hypothesis, $\gamma_i$ maps to a nonzero multiple of $f_i$ on $X_i$, so $\gamma_i\in K_{i-1}$ but $\gamma_i\notin K_i$. If $j>i$, the class $f_j$ is contracted by $X_j\to X_{j-1}$, and hence the pushforward of $\gamma_j$ to $X_i$ is zero. Thus $\gamma_j\in K_i$. Now suppose there were a nontrivial linear relation among the $\gamma_i$, and choose the smallest index $i$ with nonzero coefficient. Passing to $K_{i-1}/K_i$, all terms with larger index vanish, whereas the image of $\gamma_i$ is nonzero. This is a contradiction. The classes are therefore linearly independent. Since they all lie in $V_{X_r}^-\cap K_0$ and $\dim K_0=r$, they span $K_0$.
\end{proof}

\subsection{The reduced surface target for E5}
In the E5 model, the ordinary surface virtual dimension is negative. We therefore state the expected reduced count as a conjectural target; constructing the target family and proving compatibility of its virtual class remain open.
\begin{conjecture}\label{conj:surface-target}
Let $T=\mathbb F_2$ with negative section $S$, and let $\widetilde T$ be an ordered sequence of $N$ point blow-ups. Write $\varepsilon_0$ for the orthogonal lift of $[S]$ and $\varepsilon_j$ for the historical exceptional classes. For $\gamma=\varepsilon_0+\sum n_j\varepsilon_j$ with $\sum n_j=0$, the expected zero-dimensional virtual degree obtained by adding the transverse deformation $\mathbb F_2\rightsquigarrow\mathbb F_0$ is one if all $n_j=0$, and zero otherwise.
\end{conjecture}
\begin{remark}[Proposed mechanism and missing verifications]
Consider the projectivization of the extension $0\to\cO(-2)\to\mathscr U\to\cO\to0$ with extension class $s\eta$, where $\eta$ generates $H^1(\PP^1,\cO(-2))$. It gives a smooth projective family over $\mathbb A^1_s$ whose special fibre is $\mathbb F_2$ and whose noncentral fibres are $\mathbb F_0$. Construct the universal ordered blow-up family over this base as in \cref{lem:points}. Let $P_N$ be its parameter space and $P_{N,0}$ its central fibre. Then $P_N$ is smooth over $\mathbb A^1$ of relative dimension $2N$, while $P_{N,0}$ is irreducible. The construction also provides global markings of the historical classes by divisor classes.

Assume that a proper relative stable-map stack $M_\gamma$ with the required target-family interpretation has been constructed. For $s\ne0$, the pushforward of $[S]$ is the difference of the two ruling classes on $\mathbb F_0$, and is therefore not effective; on the central fibre, its only effective representative is $S$. The proposed argument is that a map in this class avoids the infinity section, even over Artin bases, and that its degree-one projection to $\PP^1$ survives stabilization. The key missing step is to prove that the associated affine-line-torsor section descends through stabilization and that its existence forces $s_R\eta=0$ compatibly with arbitrary Artin base change. This does not follow from the set-theoretic effectivity statement alone. If the descent and base-change claim holds, then $H^1(\PP^1_R,\cO(-2))=R$ implies $s_R=0$, so the stack factors through $P_{N,0}$.

Suppose, in addition, that the relative obstruction theory is compatible with this factorization and that the virtual class pushes forward to a cycle of the expected dimension supported on $P_{N,0}$. Write the pushforward as $\mu[P_{N,0}]$. For a general configuration, the blow-up points are distinct and disjoint from $S$. The expected geometric argument uses connectedness to rule out exceptional components; it predicts an empty fibre unless all $n_j$ vanish. In that case, the unique curve is $S$, whose normal bundle is $\cO(-2)$. The extension parameter is expected to identify with its one-dimensional obstruction space, while the remaining configuration directions are unobstructed. This predicts $\mu=1$ when all $n_j=0$ and $\mu=0$ otherwise. A proof must establish the obstruction-theory compatibility and the multiplicity, rather than deduce them from the fibrewise description alone.

Given a prescribed configuration, choose an \'{e}tale-local smooth curve in $P_N$ through that point and transverse to $P_{N,0}$. Refined Gysin base change for the relative stable-map virtual class, together with compatibility of proper pushforward and refined Gysin maps, identifies the virtual degree over this curve with its intersection number with $\mu[P_{N,0}]$, namely $\mu$. Thus the resulting degree is independent of the chosen transverse curve.
\end{remark}

\subsection{Incidence, obstruction quotient, and conditional historical preservation}
We next describe the additional steps needed to apply \cref{conj:surface-target} to arbitrary E5 historical classes. Work in the standard model $W=\operatorname{Tot}_{\PP^2}\cO(-2)$, and let $\pi:Y\to W$ be a finite sequence of point and curve blow-ups. Write $r:Y\to\PP^2$ for the projection. For the historical line class $\beta$, stabilization of the planar projection gives a morphism $\Phi:M_\beta\to B=(\PP^2)^\vee$. We may use a compactification provided the relevant maps remain in the local exceptional locus. Let $s$ be the pullback of the radial section of $\cO(-2)$. Its pullback to the universal curve vanishes because the line bundle has degree $-2$ on the unique degree-one component and degree zero on the remaining rational trees. The first direct image is $\Phi^*\cO_B(-1)$.

\begin{proposition}\label{prop:finite}
The differential of $s$ induces a surjection from the relative map obstruction sheaf to $\Phi^*\cO_B(-1)$ away from finitely many line parameters.
\end{proposition}
\begin{proof}
The union of the images in $W$ of all blow-up centres has dimension at most one, so only finitely many lines in the original plane are contained in it. For any other line, the main component meets the locus where $\pi$ is an isomorphism in a dense open subset. The differential to the radial line bundle is therefore generically surjective, and its cokernel on the main component is torsion. It follows that the induced map on $H^1$ is surjective. Normalization of the domain gives a surjection from the full $H^1$ to the $H^1$ on the main component; the target line bundle has all its $H^1$ on that component. Base change and Nakayama's lemma then give the asserted sheaf surjection.
\end{proof}
The incidence variety $I=\{(y,L):r(y)\in L\}$ is a smooth $\PP^1$-bundle over $Y$. Along a historical stable map, its relative tangent line bundle has degree $-1$ on the main component and degree zero on every other component. Its $H^0$ and $H^1$ both vanish. Consequently, lifting maps to the incidence variety does not change the relative deformation complex over the stack of prestable curves. This avoids the singular transverse surface fibres that occur at tangent lines.

For a sufficiently general line, each subsequent centre either misses the transverse surface or meets it in finitely many transverse points. The resulting surface is an iterated point blow-up of $\operatorname{Tot}\cO_{\PP^1}(-2)$. Orthogonality of the threefold historical class gives $\sum_{j\in J_i}n_j=0$ for each group of surface points arising from a threefold blow-up. Only finitely many effective surface classes can contribute: their degree against the restricted ample divisor is fixed, and their pairings with a divisor basis are bounded after adding and subtracting sufficiently large multiples of that ample divisor.

\begin{proposition}\label{prop:conditional}
Suppose that the quotient in \cref{prop:finite} descends to the absolute obstruction sheaf, that its kernel contains the intrinsic normal cone, and that the associated kernel Gysin class restricts on a general transverse surface to the actual-family reduced class predicted by \cref{conj:surface-target}. Assume also that the required properness and local-to-global comparison for the standard E5 neighbourhood hold. Then
\[
\Phi_*[M_\beta]^{\rm vir}=-H_B\cap[B],\qquad
\langle r^*H^2\rangle_{0,1,\beta}=-1.
\]
\end{proposition}
\begin{proof}
Remove the finite exceptional set from \cref{prop:finite}, and denote its complement by $B^\circ$. The assumed cone containment gives a reduced cycle of dimension two by Gysin pullback in the kernel vector-bundle stack. The excess-intersection formula yields $[M^\circ]^{\rm vir}=c_1(\Phi^*\cO_B(-1))\cap[M^\circ]^{\rm red}$. By the assumed comparison with \cref{conj:surface-target}, the reduced pushforward has generic degree one and hence equals $[B^\circ]$. The projection formula gives the required one-dimensional pushforward on $B^\circ$. Since removing finitely many points does not change $A_1(B)$, the equality extends to $B$. Finally, imposing the marked incidence condition that the map pass through a fixed point of the plane pushes forward to $H_B$, giving the numerical formula.
\end{proof}
One possible approach to the cone-containment statement is to use local trivializations of the quotient line together with the cosection lemma \cite{KL13}. The essential issue is comparison with the actual target deformation, including compatibility of the morphisms to the cotangent complexes; equality of obstruction-space dimensions alone is not sufficient. The proposition is therefore stated conditionally.

The E4 case presents a different issue. Here the obstruction is an actual quadratic equation on the ruling parameter, rather than an excess obstruction on a smooth parameter space. The local deformation $XZ-Y^2+\tau W^2+rW^3=0$ gives the ruling equation $D^2=\tau$, which separates the two rulings when $\tau\ne0$. An arbitrary curve centre need not extend in this deformation. Thus, for general sequences of curve blow-ups, one needs a relative virtual-pullback comparison over the nonreduced ruling scheme. The point-centre case is already covered by \cref{prop:historical-point}.

\end{document}